# Looking from the inside and from the outside

A. Carbone and S. Semmes [1]

June 30, 1996

## Abstract

Many times in mathematics there is a natural dichotomy between describing some object from the inside and from the outside. Imagine algebraic varieties for instance; they can be described from the outside as solution sets of polynomial equations, but one can also try to understand how it is for actual points to move around inside them, perhaps to parameterize them in some way. The concept of formal proofs has the interesting feature that it provides opportunities for both perspectives. The inner perspective has been largely overlooked.

A principal observation of this paper is that mathematical structures can be embedded into spaces of logical formulas and inherit additional structure from proofs. We shall look at finitely generated groups, rational numbers and $SL(2, \mathbf{Z})$, and examples from topology and analysis.

One of the main themes of this paper will be the possibility to embed mathematical structures into the space of logical formulas through a notion of *feasibility*. The combinatorics of proofs can then induce new structures on the original mathematical objects. We shall often see how the existence of short proofs reflects the structure of the underlying object, especially internal symmetry that is susceptible to dynamical processes.

We shall discuss these ideas for *finitely generated groups*. This will induce a new geometry, different from the one coming from the word metric. Already for free groups we shall see that the new structure is quite involved.

We shall consider *rational numbers*, where the logical notion of feasibility will be illustrated through the action of $SL(2, \mathbf{Z})$ on $\mathbf{Q}$ by projective transformations. This approach to the complexity of rational numbers suggests relationships between the lengths of proofs and number-theoretic properties.

We shall consider well-known constructions in *topology*. In spirit these constructions are compatible with the idea of feasibility. We shall describe an example of a *torus bundle* in which there is exponential distortion induced by a simple cycling which is similar to the cycling and substitutions that can occur within proofs. We shall also look at ordinary differential equations and exponentiation with continuous parameters.

[1] The first author was supported by the Lise-Meitner Stipendium # M00187-MAT (Austrian FWF) and the second author was supported by the U.S. National Science Foundation. Both authors are grateful to IHES for its hospitality.

We begin with primitive questions about the way in which mathematical objects are described, in mathematics in general and logic in particular.

# 1 Inner and outer descriptions

What is a set? How can a set be described? These are basic questions which reverberate in mathematics down to the foundations. Let us consider them here in the practical way of what mathematicians actually do.

There is a basic distinction between what one might call *inner* and *outer* descriptions of sets. For an outer description one might have a given set $A$ embedded into some larger space $X$ of simple structure, and one may describe $A$ by specifying rules which determine which elements of $X$ are in $A$. An inner description might provide a listing of the elements of $A$, with more concern for the internal structure of $A$ than an embedding of it into a larger space.

Let us consider an example. How can we describe a curve $\Gamma$ in the plane? One answer might be to provide a parameterization of it, $(x(t), y(t))$, $t \in \mathbf{R}$. Another possibility is to define $\Gamma$ as the set of solutions to some equation,

$$\Gamma = \{(x, y) \in \mathbf{R}^2 : F(x, y) = 0\},$$

where $F(x, y)$ is some function.

These are very different ways to describe a curve. In the first case it might be easy to generate many points on the curve without having a general understanding or test for when a point lies on it. For inner descriptions it may not be clear how many points are needed to have a reasonably accurate picture of the set in question, and one may have to be careful about exploring well one part while missing another. In the second case one might have a simple characterization of the elements of the set without a clear idea of how to find actual solutions.

Consider the case where we define $\Gamma$ as the zero set of $F(x, y)$ with $F$ a polynomial. A basic point about the algebraic notion of a plane curve is that it may not be compatible with the notion of a parameterization. Over the real numbers the zero set might be empty, or have several components, including compact components, etc. Some of these problems can be alleviated by working with complex numbers and making assumptions of irreducibility. A more interesting incompatibility with the idea of a parameterization is that the curve might not be rational, so that it may not be reasonable to try to parameterize the curve with the ground field. It might be an elliptic curve or a curve of higher genus. (See [29] for a discussion of algebraic curves.)

If we work over a field like the rational numbers there might be more basic problems about the existence of points on the curve. (Someone once pointed out to us that a great idea in algebraic geometry is that one can study sets of equations independently of whether one knows that there are solutions.)

Instead of thinking algebraically we can think more in terms of calculus. We should be careful about what kind of functions we allow. For instance, any closed subset of the plane can be realized as the zero set of a $C^\infty$ function $F$. (Take $F(x,y) = \exp(\Delta(x,y)^{-1})$, $\Delta(x,y)$ as on p171 of [44].) Thus the $C^\infty$ property is too flexible by itself to provide a practical way to describe sets. We can avoid this problem by restricting ourselves to smooth functions $F$ whose gradient $\nabla F$ does not vanish on the zero set of $F$. This is the hypothesis of the implicit function theorem, which then implies that the zero set of $F$ is locally given by a smooth curve. One can have singularities at critical points, as in the case of polynomials, and there is a theory for analyzing these.

This is a very basic example which hopefully illustrates well what we have in mind by "inner" and "outer" descriptions. We also see how the context matters. It is very different to think algebraically in terms of polynomials than in terms of smooth functions. Calculus permits a more flexible idea of function for which it is easier to deal with the notion of tangents naively, while algebra is more rigid but enjoys more flexibility of context, in that one can switch to the rationals or other ground fields. In algebra the singularities can be awkward to manage but their complexity is controlled at least in principle, while calculus tolerates nearly arbitrary behavior through infinite processes.

If we want to make an inner description by listing points, this listing should respect the structure of the situation, like smoothness or algebraic properties. In topology one would impose continuity conditions on mappings, in geometry and analysis one might require that mappings do not distort distances too severely.

These ideas show up in many different contexts in mathematics. For instance one can take some finite set $\mathcal{A}$ as an alphabet and look at the set $\mathcal{A}^*$ of all *words* generated by elements of $\mathcal{A}$, i.e., all finite strings of elements of $\mathcal{A}$. One might have a *language* $L$ based on $\mathcal{A}$, which is to say a subset of $\mathcal{A}^*$. A priori $L$ could be anything. How might it be described? $L$ might be *effectively enumerable*, so that there is an algorithm for generating all of its elements. This is a kind of inner description. Instead there might be an outer description, like an algorithm that says when a word lies in $L$.

Roughly speaking, inner descriptions correspond to ways to produce the effective witness, while outer descriptions correspond to ways to check mem-bership and to decide yes/no questions.

If one knows that $L$ is pretty "thick" – i.e., one has reasonably large lower bounds on the number of elements of $L$ among the words of length $n$ – then one might be able to get a practical way to list the elements of $L$ from the algorithm for deciding whether a word lies in $L$ or not. One simply goes through all the words and keeps the elements of $L$. This need not work very well if $L$ is too sparse. In the case of integer solutions of a polynomial equation $F(x, y) = 0$ it may be very difficult to tell if there are any solutions (Hilbert's tenth problem) or to know how many. In general the existence of integer solutions of polynomial equations is algorithmically undecidable, but this is not known for the rationals or other fields.

Conversely, there are sets which are effectively enumerable but for which there is no algorithm to decide membership.

There are nice variations on this theme of thickness and sparseness of languages in the context of the P = NP problem. See [28], p.87.

As another example suppose that we have an $n \times n$ matrix of complex numbers, which we think of as defining a linear mapping $T$ on $\mathbf{C}^n$. A complex number $\lambda$ is an *eigenvalue* of $T$ if $\lambda I - T$ is not invertible as a linear mapping on $\mathbf{C}^n$. The set of eigenvalues is called the *spectrum* of $T$. One can define it more concretely as the set of zeros of the polynomial equation

$$\det(\lambda I - T) = 0.$$

This is a perfectly good definition of the spectrum, but how does one actually find eigenvalues? This is a tricky question whose numerical solution is of great importance and much studied.

Dynamical systems provide another interesting case to consider. One might be able to generate a good approximation to an attractor quickly from the inside, looking at iterates of a critical point for instance, while the "rules" which govern the geometry of the attractor might be hard to see. The number of points needed to have an accurate picture of the attractor might be unclear as well.

Inner and outer descriptions need not be very compatible with each other. In mathematics one is often much more accessible than the other. Which one is more accessible depends on the context.

In the next section we shall discuss the particular case of the set of all tautologies inside the space of formulas. In Section 3 we discuss how individual formulas can in turn be used to describe subsets of other sets (equipped with some structure), and in Section 4 we consider the general relationship between algebraic structures and points inside a set.

In Section 5 we take up the notion of *feasibility*. This provides a way for us to embed a mathematical structure inside the space of formulas (with respect to some language). The combinatorics of formal proofs then induces *new* structure on our original mathematical object. For this the *cut rule* is particularly relevant. This idea is developed through examples in Sections 5-10.

This paper is intended to be accessible to a broad audience. Readers not very familiar with formal proofs may find [10] a useful source of background material.

## 2 The set of tautologies

The set of tautologies provides an interesting case to consider for inner and outer descriptions. One can consider either propositional or predicate logic.

Imagine fixing a collection of variables and the rest of a logical language, so that one has specified a notion of formulas. Let us think of the set of all formulas as being relatively simple (e.g., a recursive set), and imagine that we are interested in understanding the set of all tautologies as a subset of it through both inner and outer descriptions.

For an inner description of the space of all tautologies we can use *proofs*. The rules for building proofs provide a way to move around in the space. It may not be easy to reach a particular tautology, but in principle we can go anywhere in the space through proofs.

Given two tautologies we can always make a new one through binary logical rules. As we wander around in the space we may very well return back to the same tautology over and over again. The structure of the possible ways to move within the space reflects its geometry.

The idea of the relationship between the geometry of a space and the ability to move around in it is much studied in other parts of mathematics.

What about outer descriptions? We can use *semantics* to provide a kind of outer description of tautologies. The completeness theorem says that the set of provable formulas is the same as the set of formulas which are "true" under all interpretations. We can think of each interpretation as a test. Although there are many such tests (and indeed the set of predicate tautologies is algorithmically undecidable), it is remarkable nonetheless that tautologies enjoy these outer and inner descriptions simultaneously. One can argue that it is reasonable that neither description is very simple given that we are lucky enough to have both.

In the case of propositional logic some of these issues emerge more clearly. The "outer" characterization of tautologies as being the formulas which are true in every interpretation implies that the set of tautologies is co-NP. If P = NP (and hence P = co-NP) then there is a polynomial-time algorithm which tells whether a propositional formula is a tautology. This would be a very effective outer description.

It is not known exactly how the size of a tautology is related to the size of its shortest proof. The existence of short proofs is a way to say that the *inner* description of propositional tautologies through proofs is efficient.

Propositional and predicate logic provide very basic examples of sets in mathematics whose descriptions one would like to understand better. Another interesting example is provided by Brouwer's intuitionistic logic. In this system disjunctions and existential quantifiers are treated differently from classical logic, so that one cannot assert $A \vee B$ without actually having a proof of one of $A$ or $B$. In particular one does not take $A \vee \neg A$ as being automatic. The structure of proofs is somewhat simpler in this case than in classical logic, but the correct notion of interpretations for characterizing tautologies is more complicated.

We should mention that this is the same Brouwer who proved the famous fixed-point theorem.

Logical formulas themselves can describe mathematical objects. Our ability to describe the set of formulas is tied to the way that individual formulas can speak about underlying objects. The existence of short proofs leads to an efficient inner description of the set of all tautologies. A proof of a formula can reflect the structure of the underlying objects. Intuitionistic and classical logic differ both in the way that they describe mathematical objects and in the way that their sets of tautologies are described. For Brouwer a proof is a kind of function, where a rule like Modus Ponens corresponds to composition of functions. This is connected to the theory of Lambda Calculus, which associates functions to intuitionistic proofs in a way that reflects their internal structure.

# 3 Describing sets through logical formulas

Mathematical logic provides interesting ways to make descriptions of sets. The most basic method comes from the model theory for first-order logic. The reader who is not familiar with these concepts need not lose heart, we simply want to have an impression in mind.

In logic one has a language for talking about structures and about defining special subsets of a given set abstractly, independently of any specific set. Before we say what this means in general let us think about groups. There is an abstract idea of groups that exists independently of any particular group. There is also a way to talk about certain subsets of a group, like the elements of a certain order, or the center of the group, that exists independently of any particular group.

The idea of abstract mathematical structures can be formalized through a logical language. In a first-order language one has the usual logical connectives which represent "and", "or", "not", "implies", and the quantifiers "for all" and "there exists", but one has additional symbols which reflect the particular structure. One allows variables, constants, function symbols, and predicates (relations). Each function symbol and predicate has a fixed number of arguments, called the arity. The number and the arities of the function symbols and relations may depend on the given mathematical structure.

For example, for the theory of groups one uses one relation, the binary relation of equality $=$. There are two function symbols, one binary and one unary, which correspond to group multiplication and inversion. There is one constant symbol, corresponding to the identity element of the group.

These are all just symbols however, there is no underlying set. That is because a first-order language intends to describe the *idea* of a group rather than a particular one.

One also needs the notion of a *term*, which is an expression constructed from variables and constants using function symbols. Think of a formal expression for groups, some product of variables, possibly with inverses. The functions and relations take terms for their arguments, as in the composition $s^{-1}t$ and the relation $s = t$ in the context of groups, where $s$ and $t$ are terms.

A relation with a choice of arguments – like $s = t$ – is a logical formula in a first-order language, an atomic formula. Informally it is a statement which might be true or false, depending on the context. These atomic formulas can be combined with the logical connectives to build more complicated formulas.

The *theory* of groups consists of the usual axioms governing the group operations, such as the associativity axiom.

All of this exists purely at the level of formal symbols. Roughly speaking a *model* is a specific interpretation of the constants, functions, and relations over some set. This means a set $S$ in which the variables can take values, with particular choices in the set for constant symbols and particular choices for the functions and relations. These choices should satisfy the axioms of the theory. Thus actual groups are models for the first-order theory of groups, with an actual set of group elements, the usual notion of equality, a choice of group operations, etc. There are many different kinds of groups, many different models, but just one first-order theory of groups.

Another example is provided by arithmetic. One can formalize it with the binary relations $=$ and $<$, operations like addition and multiplication, and the well-known Peano axioms. The usual notion of natural numbers provides a model for this theory, but there are nonstandard models too.

A first-order language provides the possibility to make universal recipes for describing certain sets, a set for each choice of model. A formula in the language, like $\phi(x, y, z)$, with the free variables $x$, $y$, and $z$ and no others, defines such a recipe. Given a specific model based on a set $S$ we can get a subset of $S \times S \times S$, namely the set of triples $(x, y, z)$ for which $\phi(x, y, z)$ is a valid formula. For instance one can define the center of the group in this way, or the set of elements of order 2. In the theory of fields one can define algebraic varieties without specifying the particular ground field. For ordered fields one can define sets in terms of polynomial equations and inequalities.

One can think of this as a way to make "outer" descriptions of certain classes of sets through logic. It is rather sophisticated, because of the possibility of quantifiers. Without quantifiers it is already subtle, a kind of tricky extension of algebra, but quantifiers make it even more complicated. Indeed, many problems of decidability involve finding a uniform way to eliminate quantifiers.

These are well-known and much studied ideas, particularly in connection with algebraic sets and their generalizations. However the model-theoretic view does not seem to provide good ways to explore sets from the inside.

## 4  Some comments about algebra and points

The preceding discussion is reminiscent of a general phenomenon in algebra, in which one has algebraic structures which make sense abstractly but which arise classically from actual sets with actual points. A basic example of this is given by the following. Let $X$ be a compact Hausdorff topological space, and let $C(X)$ denote the space of all complex-valued continuous functions on $X$. This is an algebra, and even a $C^*$-algebra. However one can talk about algebras abstractly even if they do not arise from actual spaces in this way.

We would like to say that having a theory (in the sense of logic) is like having an algebra without necessarily having a space of points underneath. One can have *representations* of the algebra on spaces with actual points, as in the notion of a structure in model theory.

The idea of "inner descriptions" should entail actual points in actual spaces and their interactions. To what extent can we do that directly from the algebraic structure? For this question the example of abstract algebras and algebras of functions is instructive.

If $X$ and $C(X)$ are as above, then we can recover the points in $X$ from the algebraic structure in $C(X)$. If $p$ is a point in $X$, then $\{f \in C(X) : f(p) = 0\}$ is a maximal ideal in $C(X)$, and conversely every maximal ideal arises in this manner. Given another compact Hausdorff space $Y$ one can show that $X$ is homeomorphic to $Y$ if and only if $C(X)$ is isomorphic to $C(Y)$ as an algebra.

One can characterize the algebras that arise this way, as commutative $C^*$-algebras. See [38, 40]. There are analogous stories in the context of algebraic varieties, but let us stick to topological spaces for simplicity.

In principle compact Hausdorff spaces are described completely by the algebra of commutative functions on them, but how does this work practically? How can one see inside the space through the algebra? In some kind of practical way, and not just in principle? This turns out to be very mysterious, and not much is known. There is a different way to try to represent the structure of a space in purely algebraic terms, through which one can recover topological invariants of the underlying space from direct algebraic constructions. See [13]. This approach also gives meaning to these topological invariants in non-commutative settings where there need not be "points" in the classical sense, and this is a matter of great current interest in mathematical physics.

This is similar in spirit to the relationship between operational and denotational semantics in programming languages. See [47].

It can happen naturally that one has an algebra in hand but not the underlying points that one wants. For example, let $T$ be a linear transformation acting on some $\mathbf{C}^n$. Consider the algebra of linear transformations generated by $T$, which amounts to saying all polynomials in $T$. This is a nice commutative algebra, but what are the underlying "points"?

Suppose that $T$ is diagonal and has distinct eigenvalues $\lambda_1, \ldots, \lambda_n \in \mathbf{C}$. These eigenvalues are the "points" in a natural way. If we let $X$ denote the set of them, then the algebra of linear transformations generated by $T$ is isomorphic to the algebra of polynomials restricted to the set $X = \{\lambda_1, \ldots, \lambda_n\} \subseteq \mathbf{C}$ of eigenvalues.

If $T$ is not diagonal but is diagonalizable then the algebra which it generates has the same form. If it is not diagonalizable, having a nontrivial Jordan canonical form, then the notion of "points" underlying the algebra is more problematic, because of nilpotent elements.

There are versions of this discussion for linear operators acting on infinite-dimensional spaces, in which the natural notion of spectrum is an infinite set whose topological structure becomes important. See [38]. This brings us closer to the earlier discussion.

A better analogy between logical theories and algebraic structures in which the concept of "points" is not directly involved is provided by Boolean algebras. In the Boolean case the elements of the algebra are closer to points in an underlying set (in the sense of Stone's theorem) than for algebras of complex-valued functions.

In proofs the algebraic structure is coded through rules of inference. It is also natural to think of proofs simultaneously as being sets with points (such as atomic formulas) and combinatorial structure. This is a remarkable coexistence.

This idea is illustrated in a strong way by [6], in which the *Craig interpolation theorem* [15] is discussed in a combinatorial context without the algebraic structure of connectives. This combinatorial view of points in a proof is also present in the notion of *logical flow graphs* [5, 7], which trace the logical connections within a proof, and in the study of *cycles* in these graphs, as in [7, 8, 9]. (Logical flow graphs are related to the earlier notion of *proof nets* [20].) The concept of *inner proofs* from [7] reflects the idea that there are points walking around inside proofs.

We shall give here other examples to illustrate this idea of proofs dealing with physical points and not just algebraic constructions. One of the tools for doing this will be the notion of *feasibility*, which arose from an extension of arithmetic that we shall discuss next.

## 5 Feasible numbers

There has been much concern in mathematics about abstraction which may not reflect anything concrete or "real". Extremely large numbers were troubling to some, and there was the idea that they should be treated differently from a small number like 37 which is closer to ordinary existence.

The first mathematical treatment of *feasible numbers* was given in [36]. (The philosophical discussions go back to Mannoury, Poincaré, and Wittgenstein.) For this we start with the first-order theory of arithmetic, and we add a unary predicate $F$. Roughly speaking $F(x)$ is interpreted as meaning that $x$ can be constructed in some feasible manner. We shall use the arithmetic operations $+$, $*$ (multiplication), and $s$ (successor). In addition to the usual axioms of arithmetic we add the following axioms for $F$:

$$
\begin{array}{ll}
& F(0) \\
F : equality & x = y \rightarrow (F(x) \rightarrow F(y)) \\
F : successor & F(x) \rightarrow F(s(x)) \\
F : plus & F(x) \wedge F(y) \rightarrow F(x + y) \\
F : times & F(x) \wedge F(y) \rightarrow F(x * y)
\end{array}
$$

In other words, 0 is considered to be feasible, and the property of feasibility is closed under equality, successor, addition, and multiplication.

For this discussion we do not permit ourselves to use induction over $F$-formulas. Otherwise we could prove $\forall x F(x)$ in a few steps. Note that if we add the axiom $\exists x \neg F(x)$, asserting the existence of a nonfeasible number, then we still get a consistent system, for which the models are nonstandard models of arithmetic.

The idea instead is that if we can write down a proof of $F(t)$ for some term $t$, then that should mean that $t$ was "feasible" in a reasonable sense. Of course we can always prove $F(n)$ for any natural number $n$ in about $n$ steps, using the successor rule repeatedly. (Strictly speaking we are abusing the first-order language of arithmetic here, $n$ really means the result of applying $n$ times the successor function to 0. Syntactic technicalities only detract from the main points and we shall ignore them without remorse.) However we can use the size of a proof of $F(n)$ as a measurement of the feasibility of $n$.

This is a very nice point. We can use proofs to make descriptions of mathematical objects and to make measurements of their complexity. We shall push aside the foundational issues and simply use the idea of feasibility as a tool for studying mathematical structures.

To make precise the measurements one should be careful about the formalization of proofs. We shall not discuss this in detail, but there are a couple of important points. The first is that we consider only proofs in which the result of any intermediate step is used only once. Thus proofs have tree-like structures. The second concerns the role of the "cut" rule in sequent calculus and its counterpart in other systems. Roughly speaking, the cut rule allows indirect reasoning through lemmas. It is a generalization of the deduction rule Modus Ponens, which says that if you know $A$ and if you know that $A$ implies $B$ then you can conclude $B$. Without the cut rule a proof of $F(t)$ for some term $t$ would have to exhibit an explicit construction of the term $t$. With the cut rule one can make short proofs of feasibility which provide only implicit descriptions, as we shall soon see. There are effective methods for converting proofs with cuts into proofs without, at the cost of great expansion in the proofs. See [21, 45, 10].

Let us mention one more point. In [36] an $F : inequality$ rule is included in the axioms, to the effect that if $y$ is feasible and $x < y$ then $x$ is also feasible. For the historical concern about large numbers this is a reasonable requirement to consider, but we have omitted it intentionally. It does not fit as well with the idea of a proof of feasibility of $F(t)$ as providing a description of $t$, and it is less convenient for other mathematical contexts. So we simply drop it. This is an important conceptual point. In mathematics we can make definitions to suit our purposes, irrespective of historical traditions.

So let us now consider the concrete matter of how we might give short proofs of feasibility of numbers. We follow the examples in [8].

We can always get a proof of $F(n)$ in about $n$ lines through repeated use of the $F : successor$ rule. We can be a little more intelligent and get a proof of $F(2^n)$ in about $n$ lines using the $F : times$ rule repeatedly. This is the standard mathematician's trick of using geometric progressions.

We can improve on this as follows. We know that

(1) $$F(x) \to F(x^2)$$

In particular we have that

(2) $$F(2^{2^j}) \to F(2^{2^{j+1}})$$

for all $j = 0, 1, 2, \ldots$ We can combine $n-1$ copies of (2) together with the feasibility of 2 (i.e., $F(s(s(0)))$) to get a proof of $F(2^{2^n})$ in $O(n)$ lines.

In this argument we won an exponential over the previous one. The price for this was that we were implicitly using cuts and contractions to make the building blocks and to combine them. The proof of feasibility did not furnish a direct construction. In terms of sequent calculus one can see the importance of the contractions. To prove $F(x) \to F(x^2)$ one contracts two copies of $F(x)$ on the left-hand side of the sequent into one.

We can win another exponential using quantifiers. We can prove the feasibility of $2^{2^{2^n}}$ in $O(n)$ lines. The proof is constructed from the following building blocks. First we have that 2 is feasible, as above. Next we have that

(3) $$\forall x (F(x) \to F(x^2)),$$

which was the main ingredient in the preceding construction. The last building block is

(4) $$\forall x (F(x) \to F(x^k)) \to \forall x (F(x) \to F(x^{k^2}))$$

That is, we can use $\forall x (F(x) \to F(x^k))$ twice, the second time replacing $x$ with $x^k$, to get $\forall x (F(x) \to F(x^{k^2}))$. This is much better than before, we are squaring the exponent instead of multiplying it by 2. By combining a series of these last building blocks, with $k = 2^{2^j}$, $j = 0, 1, \ldots, n-1$, we can conclude that

(5) $$\forall x (F(x) \to F(x^2)) \to \forall x (F(x) \to F(x^{2^{2^n}})),$$

and then we combine with the other pieces to get a proof of $F(2^{2^{2^n}})$ in $O(n)$ lines.

This last approach has some interesting features. As observed in [8], we only used the $F : times$ rule once, in the proof of (3). In (4) we made progress by making a substitution. The proof of (4) uses contraction rules in an interesting way, one uses $\forall x (F(x) \to F(x^k))$ twice to get $F(x) \to F(x^k)$ and $F(x^k) \to F((x^k)^k)$, but the quantifier rules permit us to convert this into two copies of $\forall x(F(x) \rightarrow F(x^k))$ which can then be contracted into each other. This is a standard point about quantifiers and contractions, they permit us to contract two occurrences of a formula into one even though they have very different histories within the proof. This kind of substitution did not occur in the previous propositional argument.

Notice that in this last proof we did not have nesting of quantifiers. There are more elaborate proofs, due to Solovay, which use many nested quantifiers to get short proofs of very large numbers defined through towers of exponentials of arbitrary height. One gains an extra exponential with each nested quantifier. See [8] for details and a discussion of the dynamical structure of these proofs.

With these examples in mind let us think about the type of description of a number provided by a proof of feasibility. In a proof of $F(2^n)$ in about $n$ lines using multiplications we really make an explicit construction. We cannot expect to do better than win an exponential, because our most powerful operation is multiplication.

The other arguments are increasingly less explicit, because of the use of substitutions. There is a kind of balance in this; as the proofs become shorter, their internal structure becomes correspondingly more complicated, reflecting the increasing difficulty by which the implicit descriptions can be unwound into explicit constructions. The *logical flow graphs* of the proofs, which trace the flow of occurrences of formulas in proofs, become more complicated, with increasing numbers of bridges and cycles. See [8, 9].

The procedure of *Cut-elimination* (introduced by Gentzen [18, 19], see [21, 45, 10]) furnishes a general method for transforming implicit descriptions into explicit constructions with effective bounds.

In giving short proofs of feasibility of large numbers like $2^{2^n}$ or $2^{2^{2^n}}$ we are using the special structure of these numbers, a kind of internal symmetry to them. This internal symmetry is reflected in the existence of short proofs, but there are no theorems about this. In general we should not be able to win so much compression using cuts, because arbitrary numbers will not have so much internal symmetry. General quantitative results have not been given.

The mathematical idea of feasibility provides a way to embed arithmetic inside a space of formulas. Formal proofs then lead to new structure for natural numbers. This structure is quite different from the ones that are usually considered, and the cut rule plays an important role in this.

# 6 Groups

In recent years much attention has been devoted to the study of the structure in finitely generated groups which can be seen through the word metric. (See [26], for instance.) One fixes a generating set and defines the distance from an element $g$ of the group $G$ to the identity $e$ to be the minimal length of the word that represents $g$. This can be extended to a left-invariant metric on all of $G$.

We can try to make other kinds of uniform measurements in the theory of groups using proofs and the idea of feasibility. Again let us introduce a unary predicate $F$, acting now on elements of our given group $G$. Let us also fix a finite subset $S$ of $G$ – we can think of it as a generating set, but actually the concept makes sense in any case – and require that $F$ have the following properties:

$$
\begin{aligned}
& F(e) \\
& F(\gamma) \quad \text{for each } \gamma \in S \\
F : equality \quad & x = y \to (F(x) \to F(y)) \\
F : composition \quad & F(x) \wedge F(y) \to F(xy) \\
F : inverse \quad & F(x) \to F(x^{-1})
\end{aligned}
$$

Here we write $xy$ for the group composition and $x^{-1}$ for the group inverse.

The length of the shortest proof of the feasibility of an element of $G$ can be taken to be some kind of measurement of its complexity. It is a well-defined function on $G$ because of the $F : equality$ rule. The length is always bounded by a constant multiple of the distance to $e$ in the word metric. We can make examples of proofs of feasibility which parallel the ones in the previous section. If $x \in G$ is feasible, then we can make a proof of $F(x^n)$ in $O(n)$ lines, by repeated use of the $F : composition$ rule. We can be more clever and get a proof of $F(x^{2^n})$ in $O(n)$ lines by making proofs of

$$(6) \qquad F(x^{2^j}) \to F(x^{2^{j+1}})$$

for $j = 0, 1, 2, \ldots, n-1$ as in (2) and combining them. This argument requires only propositional logical rules. If we use also quantifiers then we can get a proof of $F(x^{2^{2^n}})$ in $O(n)$ lines as before. That is, the proof that we outlined before for (4) works just as well here.

The last method that we mentioned in Section 5, based on nesting of quantifiers, does not work in the theory of groups. To apply it to get a universal nonelementary distortion in groups (i.e, short proofs of $F(b) \to F(b^N)$ with $N$ a tower of exponentials like $2^{2^{2^{2^2}}}$) we would have to permit ourselves to quantify over integers as well as group elements. Indeed, for this argument we need to make substitutions into exponents, and this means substitutions with integers. The other arguments require only substitutions of group elements.

At any rate we can make short proofs of feasibility using the first two methods. Given a finitely generated group these methods can be combined with the cancellation induced by the relations to yield even shorter proofs of feasibility. For example, let $G$ be the group with generators $x$ and $y$ and the single relation $y^2 = xyx^{-1}$. Thus $y^{2^m} = x^m y x^{-m}$. The feasibility of $x^m$ implies that of $x^m y x^{-m}$, and combining this with the earlier arguments we can get a proof of the feasibility of the group element $y^{2^{2^{2^n}}}$ in $O(n)$ lines. (See [26] for other examples of finitely presented groups with distortion.)

Although in a sense we are simply transferring the earlier arguments for integers (from Section 5) to the theory of groups, there is an important difference between the two situations. In groups there are many ways to go to infinity. In a free group, for instance, every infinite word describes a path to infinity in the associated tree. Our arguments about the integers lead to a lot of compression for proofs of feasibility along the direction of a cyclic subgroup $\{a^n\}$, at least for some $n$'s. In the word metric all directions towards infinity in the free group are practically the same, but in the geometry of feasibility the cyclic subgroups are very special compared to generic directions. One can think of feasibility as providing a way to measure the amount of algebraic structure in a given direction.

This point can be seen in broader terms. The amount of compression that one can get for a notion of feasibility in some context can be seen as a measurement of the internal structure of the object in question. The examples in Section 5 reflect the internal symmetry in the case of arithmetic. We can form $z^x$, where $z$ and $x$ are both integers, and we can make substitutions between them.

In the spirit of automatic groups (see [16]) one can be interested in representing a group through its set of words. We can enhance the notion of feasibility to be sensitive to the different ways that a group element is represented by words. Suppose now that our group $G$ is finitely presented, with a finite set $R$ of relations $w_i = e$ which express the triviality of the words $w_i$. We can introduce a new unary predicate $T$ so that $T(w)$ is intended as meaning that $w$ represents a trivial word. We impose the following axioms:

$$\begin{array}{ll} & T(e) \\ & T(w_i) \quad \text{for each } w_i \in R \\ T: equality & w = u \to (T(w) \to T(u)) \\ T: composition & T(w) \wedge T(u) \to T(wu) \\ T: inverse & T(w) \to T(w^{-1}) \\ T: conjugation & T(w) \to T(vwv^{-1}) \end{array}$$

The idea of the last rule is that a trivial word conjugated by any word should again be trivial. It seems reasonable to make this rule without requiring that $v$ be feasible, but one might want to make different choices, depending on the context.

For the purposes of making measurements in groups, one can combine the axioms for $F$ and $T$ and also add

$$F(w), T(u) \to F(wu) \quad \text{and} \quad F(w), T(u) \to F(uw).$$

Now $F$ is defined over words instead of group elements.

The notions presented in this section are not necessarily canonical or fixed. For instance, one might want to study chains of subgroups, each normal in the larger one, with different predicates for the different subgroups, each predicate axiomatized as above.

The bottom line is that proofs provide a nice way to try to look inside groups, to move around inside them and test their structure. This is an idea that has not been explored.

The view of groups and proofs described in this section was motivated by [9], which goes in the opposite direction: one starts with a proof and associates a group to it to reflect its structure.

One feature of the idea of feasibility is its universality. It applies to all groups at once, like the word metric. This universality continues to exist under restrictions on the kind of proofs that we allow. This is an important point: one is free to choose a fragment of logic to suit one's purposes. Different fragments can lead to different metrics on groups.

# 7  Rational numbers

We can extend the idea of feasibility to rational numbers in a natural way. For our purposes it will be convenient to consider $\infty$ as a rational number, with the conventions that $\infty \cdot \infty = \infty$, $a \cdot \infty = \infty \cdot a = \infty$ when $a \neq 0$, $0 \cdot \infty = \infty \cdot 0 = 0$, and $\frac{a}{\infty} = 0$ when $a$ is a finite rational number. We leave all other cases undefined. The need for $\infty$ is slightly a nuisance, but the point of it will be clear in a moment, and these technicalities are not serious.

We can introduce a feasibility predicate $F$ in much the same way as before. Now we want to measure "rational" complexity, and we want to take the field structure into account. We give ourselves the obvious rules for $F$, namely that 0 and 1 are feasible, that equality, sums, and products preserve feasibility, and that additive and multiplicative inverses preserve feasibility. There are some small caveats needed to account for the cases when the operations are not defined, but let us not worry about that.

We have seen how the feasibility of large integers can be established through short proofs, we can do the same now for rational numbers. This suggests a very natural open problem: how can one relate the number-theoretic properties of a rational number to the size of proofs of feasibility? Continued fraction expansions might be natural for this. Notice that we can get short proofs of the feasibility of $\frac{2^m}{3^n}$ and $\frac{2^m+5^j-12}{3^n+7^k}$ for large $j, k, m$, and $n$ as in Section 5. For proving feasibility the two are practically the same, but for number theory they are quite different.

As for groups, the restriction to different fragments of logic can lead to different number-theoretic properties, and one is free to choose the logical system to suit one's purposes.

Let us describe now an amusing construction for feasibility of rational numbers. The basic point is that $2 \times 2$ matrices with (finite) rational entries act on rational numbers in a natural way. Let $\begin{pmatrix} a & b \\ c & d \end{pmatrix}$ be such a matrix, and consider the transformation

$$(7) \qquad x \mapsto \frac{ax+b}{cx+d}$$

We assume that the determinant of our matrix is different from zero to avoid problems with the definition. This condition ensures that the numerator and the denominator above cannot both vanish at the same time, so that the quotient is always defined. It is for this reason that we allow $\infty$ as a rational number. If $x = \infty$ then we interpret the above quotient as being $\frac{a}{c}$. Not both of $a$ and $c$ can vanish, because of the assumption of nonzero determinant.

Let $A$ denote such a matrix $\begin{pmatrix} a & b \\ c & d \end{pmatrix}$, and let $A$ also denote the projective linear transformation defined in (7). The correspondence from matrices to projective transformations is a homomorphism, as is well known. Indeed, in working with the rational numbers $\mathbf{Q}$ together with $\infty$ we are really working with the projective line over the rational numbers, which means the space of ordinary lines in $\mathbf{Q} \times \mathbf{Q}$. If $\alpha \in \mathbf{Q}$, then we associate to it the line that passes through $(1, \alpha)$. This parameterizes all lines in $\mathbf{Q} \times \mathbf{Q}$ except for the one which passes through $(0, 1)$, which we associate to $\infty$. A matrix $A$ acts linearly on $\mathbf{Q} \times \mathbf{Q}$ and induces a transformation on the space of lines in $\mathbf{Q} \times \mathbf{Q}$. This works out to be compatible with the mapping above.

Let $x$ be a rational number and consider $A^n x$. We would like to have short proofs of feasibility of $A^n x$ for large values of $n$. This fits with the earlier construction for groups. Here we want to use instead the notion of feasibility for the field of rationals. We encode feasibility for the group of matrices into the feasibility for the field as follows. Given a $2 \times 2$ matrix $A$ with finite rational entries, let us write $\phi(A)$ for the formula $F(a) \wedge F(b) \wedge F(c) \wedge F(d)$. This extension of feasibility is preserved by *matrix* multiplication, by an easy argument.

This permits us to make short proofs of $\phi(A^n)$ for large values of $n$ as in Section 6. For the argument in Section 5 based on nested quantifiers there are some subtleties. To understand the issue properly we should first observe that the constructions discussed here can be extended to fields in general. To make the argument using nested quantifiers we need to be able to quantify over integers (which would arise in the exponents of our matrix). One can do this if one can *define* the integers inside the field. For the rationals there is a way to do this, due to Julia Robinson [37]. This would not work in a field of finite characteristic.

Once we have short proofs of $\phi(A^n)$ for large $n$ we can get short proofs of the feasibility of $A^n x$ for large $n$, given the feasibility of $x$.

We chose this example in part because of the well-known importance of projective transformations in analysis and number theory. Let us review some aspects of complex analysis and its connection with rational numbers. In analysis one works with complex numbers, both as matrix entries and for the domain on which the projective transformations act. Instead of having them act on the whole complex plane one often restricts oneself to actions on the upper half-plane

$$\{z \in \mathbf{C} : z = x + iy, \ x, y \in \mathbf{R}, y > 0\}.$$

A well known corollary of the uniformization theorem [1, 2] in complex analysis implies that most Riemann surfaces can be realized as the quotient of the upper half-plane by a discrete group of projective linear transformations. "Most" means all Riemann surfaces except the sphere, the plane, the plane with one puncture, and tori.

These group actions on the upper half-plane induce group actions on the boundary, the real line, which should be completed by the addition of a point at infinity. The action on the boundary can be much more chaotic than in the interior, with the orbit of a point being dense instead of discrete.

Sometimes Riemann surfaces and the corresponding groups of projective transformations have additional arithmetical structures. It can be natural to look at the action of the groups on rational numbers.

Here is a special case which provides an important example. Let $SL(2, \mathbf{Z})$ denote the group of $2 \times 2$ matrices with integer entries and determinant 1. This is indeed a group. The main point is that the inverses of such matrices still have integer entries, because of the determinant. In fact it is well known to be finitely generated with generators $\begin{pmatrix} 0 & -1 \\ 1 & 0 \end{pmatrix}$ and $\begin{pmatrix} 1 & 1 \\ 0 & 1 \end{pmatrix}$. (See [31], p30, Theorem 1.) This group acts on the upper half-plane by projective transformations, and the quotient that results is isomorphic as a Riemann surface to the twice-punctured plane $\mathbf{C} \backslash \{0, 1\}$.

For the purpose of making rational numbers which admit short proofs of their feasibility we can extend the preceding construction. Because the correspondence between matrices and projective transformations is a homomorphism, powers of projective transformations correspond to powers of matrices, and we can work directly with square matrices of any rank. Suppose that $A$ is an $m \times m$ matrix with rational entries. As before we can define a formula $\phi(A)$ which expresses the feasibility of the entries of $A$, and there is a simple theorem to the effect that $A \cdot B$ has feasible entries if $A$ and $B$ do. This implies that we can make short proofs of the feasibility of matrices $A^n$ for large $n$. The entries of $A^n$ are rational numbers which have short proofs of feasibility. Note however that the method of projective transformations enjoys more structure than the case of matrices of arbitrary rank, defining a group action on $\mathbf{Q}$ in particular.

We have described now some fast constructions of rational numbers, but one can imagine plenty of others. Proofs provide a context for certain kinds of dynamical behavior. Our extensions of the constructions in Section 5 work whenever we have an action of the semigroup of nonnegative integers. The precise nature of the dynamical behavior is not yet well understood, nor is it properly reconciled with current knowledge of dynamical systems.

# 8  A story from topology

The examples so far illustrate how one might use the idea of feasibility to make measurements and descriptions of mathematical constructions through proofs. The examples all had a kind of discreteness to them, we would like to be more ambitious now and point towards a more "continuous" setting. For this we shall need to be even more relaxed than usual about formalization. Certain ideas will emerge even if nothing is precise about the logic.

The concept of feasibility has a certain affinity for continuity. One could say that it wants to provide a kind of quantitative version of connectedness.

Our example from topology will take some time to explain and so we describe some general points first. We shall begin by reviewing the concept of *Serre fibrations* from topology [39, 3]. This notion entails the construction of continuous families of mappings for which the idea of feasibility is potentially relevant. To bring out this point we shall discuss in some detail a particular example of a *torus bundle*. In this special case the required construction amounts to taking large powers of a matrix in $SL(2, \mathbf{Z})$. The idea of feasibility seems to be compatible with Serre fibrations in more general situations, though, even if it does not reduce to groups as easily as in the example of torus bundles over the circle. In the smooth case, for instance, one can make constructions by solving ordinary differential equations, which can be seen as a generalization of taking large powers of a matrix.

Our example of a torus bundle captures geometrically a basic phenomenon in proofs. Sometimes complicated constructions can be coded in short proofs through repeated cycling and substitutions. A proof may describe a simple operation which is used repeatedly in the actual construction. In our topological example we shall see that the simple motion of cycling around a circle many times induces a motion up in our torus bundle with exponential distortion.

Let us now proceed with the details. Let $E$ and $B$ be two topological spaces, and let $\pi : E \to B$ be a mapping between them. We say that $\pi : E \to B$ is a *Serre fibration* if it enjoys the following property. Let $P$ be a finite polyhedron, and suppose that we have continuous mappings $f : P \to E$ and $g : P \times [0,1] \to B$ such that $\pi \circ f = g(\cdot, 0)$. Then there should be a "lifting" $\widehat{g} : P \times [0,1] \to E$ of $g$, meaning a continuous mapping with $g = \pi \circ \widehat{g}$, such that $\widehat{g}(\cdot, 0) = f$. This is similar to the lifting of paths in covering surfaces [2, 33], but now we are working with continuous families of paths parameterized by the polyhedron $P$.

To understand what this means consider the simple case where $E = B \times F$ for some topological space $F$. (Here $B$ is the "base" and $F$ is the "fiber".) In this case the fibration property is automatic, one can write down a $\widehat{g}$ directly. (Take $\widehat{g}(p, x) = (g(p, x), \pi_1(f(p)))$, where $\pi_1 : E \to F$ is the obvious projection onto $F$.) In general, topologists are more interested in situations where $E$ looks like a product above small subsets of the base $B$, but for which there is nontrivial twisting around globally. In these cases one can often still verify the existence of the necessary lifting. The point is to exhaust $P \times [0, 1]$ through local liftings. For this one needs compactness assumptions to ensure the finiteness of the construction.

Serre fibrations are useful because one can relate the topology of the total space $E$ to the topology of the base $B$ and the fibers through an exact sequence of homotopy groups. This is most interesting in the case where $E$ is not simply a product, so that there is nontrivial twisting.

Let us think of the lifting $\widehat{g} : P \times [0,1] \to E$ of $g$ above as being like an explicit proof of feasibility. We start with an initial configuration which is given by $f : P \to E$, and $\widehat{g}$ provides a way to get from $f$ to a final configuration given by $\widehat{g}(\cdot, 1)$ in a continuous manner. In a discrete context one would think of a sequence of very small steps. In practice there can be a definite amount of distortion at each step, and then exponential distortion over the whole time interval $[0, 1]$. See the figure below.

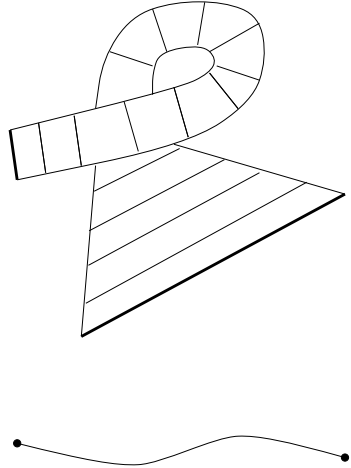

One can be interested in knowing how hard it is to find this lifting, how complicated the proof has to be for instance. This is not at all well defined because we have not been precise about the formalization of logic in this case, and it is not exactly clear what operations should be allowed. As a matter of principle it does fit well with our earlier examples, though, because there are situations where the lifting is obtained by iterating a simple operation. It is easier to think about this in the setting of compact smooth manifolds, in which the construction can be made by solving an ordinary differential equation, which is a continuous version of iteration. (For this a technical assumption on the mapping $\pi$ is needed, which is that it be a *submersion*, which means that its differential is surjective everywhere.)

To make the ideas and the role of feasibility more clear let us consider a concrete example of a fibration with nontrivial twisting. Let $S^1$ denote the unit circle. It will be more convenient to think of it as the quotient space $\mathbf{R}/\mathbf{Z}$. Let $T$ denote the torus $S^1 \times S^1$, which we can think of as $\mathbf{R}^2/\mathbf{Z}^2$.

We want to look at *torus bundles* over a circle. We shall use the following recipe. Suppose that $A : T \to T$ is a homeomorphism. Take $[0, 1] \times T$ and glue the two ends $\{0\} \times T$ and $\{1\} \times T$ together using $A$. This means that we take $[0, 1] \times T$ and we identify $(0, u)$ with $(1, A(u))$ for all $u \in T$. This defines a space which we call $E$. There is a natural mapping $\pi : E \to S^1$ which corresponds to the projection of $[0,1] \times T$ onto $[0,1]$, where we identify $S^1$ with the space obtained by taking $[0,1]$ and identifying the endpoints $0$ and $1$.

We can describe this space in another way as follows. We start with $\widetilde{E} = \mathbf{R} \times T$. We define a mapping $\phi : \widetilde{E} \to \widetilde{E}$ by $\phi(x,u) = (x+1, A(u))$. This mapping generates an infinite cyclic group of homeomorphisms on $\widetilde{E}$, and $E$ is just the quotient of $\widetilde{E}$ by this group, in the same way that $S^1$ is the quotient of $\mathbf{R}$ by the infinite cyclic group of homeomorphisms generated by $x \mapsto x+1$.

If instead of $\phi$ we used the mapping $(x,u) \mapsto (x+1, u)$ we would simply get $S^1 \times T$ for the quotient. By choosing a suitable mapping $A : T \to T$ we can get a bundle in which there is some nontrivial twisting as we go around the base.

Let us consider now a specific example of such a mapping $A$. Start with the matrix $\begin{pmatrix} 2 & 1 \\ 1 & 1 \end{pmatrix}$, which lies in $SL(2, \mathbf{Z})$ since it has determinant 1. This defines a linear mapping on $\mathbf{R}^2$. Because the matrix has integer entries the corresponding linear mapping sends the standard integer lattice $\mathbf{Z}^2$ inside $\mathbf{R}^2$ to itself. Thus we can get a well-defined mapping on the quotient $\mathbf{R}^2/\mathbf{Z}^2 = T$, and we take this to be $A$. This defines a homeomorphism on $T$, because the inverse of $\begin{pmatrix} 2 & 1 \\ 1 & 1 \end{pmatrix}$ is also a matrix with integer entries (since the determinant is 1), and hence it descends to a mapping on $T$ as well.

Thus we get a homeomorphism $A : T \to T$. It may seem harmless but in fact it is quite nontrivial. It is not homotopic to the identity, for instance. For if it were, its lifting to the universal covering of $T$ would differ from the identity by only a bounded amount, and this is not true. Indeed, $\mathbf{R}^2$ is the universal covering of $T$, and the lifting of $A$ to it is just the linear transformation that we started with.

The first homology and homotopy groups of $T$ are isomorphic to each other and to $\mathbf{Z}^2$. By general nonsense, the homeomorphism $A$ induces an automorphism on this group, which in this case is given by the action of the matrix $\begin{pmatrix} 2 & 1 \\ 1 & 1 \end{pmatrix}$ with which we began. This provides a topological way to measure the difference between $A$ and the identity mapping, even up to homotopy.

To understand better the nontrivial effect of $A$ it is helpful to compute the eigenvalues of our matrix. These are the roots of the polynomial

$$\det(\begin{pmatrix} 2 & 1 \\ 1 & 1 \end{pmatrix} - \lambda I) = (2-\lambda)(1-\lambda) - 1 = (\lambda - \frac{3}{2})^2 - \frac{5}{4},$$

namely $\lambda = \frac{3 \pm \sqrt{5}}{2}$. Note that the product of these numbers is 1, as it should be, and one is larger than 1 and the other is smaller than 1. In fact the larger eigenvalue is between 2 and 3.

Because our matrix is symmetric we can find an orthogonal basis with respect to which it is diagonal with these two eigenvalues. When we take large powers of the matrix we get exponential compression in one direction and exponential expansion in the other. Topologically this exponential expansion for the matrix implies that we can find loops in the torus $T$ whose image in $T$ under $A^n$ wraps around an exponentially larger number of times. This wrapping is depicted by the diagonal lines in the figure below, where the torus is obtained from the square by identifying the opposite sides. The diagonal lines represent a single curve in the torus.

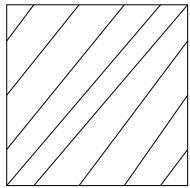

Thus $A$ is quite nontrivial and this leads to nontrivial twisting of the fibration $E$. We want to see how this twisting is reflected in liftings (as for Serre fibrations). Suppose that $f : P \to E$ is some continuous mapping, where $P$ is a finite polyhedron. One can think of $P$ as being a polygonal circle. Let us assume that the image of $f$ lies in a single fiber of $E$, so that there is a point $b \in S^1$ such that the image of $f$ lies in $\pi^{-1}(b)$. Let us assume also that our mapping $g$ is of a particularly simple form, that it is constant on $P$. Thus $g$ is in essence a mapping from $[0,1]$ into $S^1$, which we denote by $\gamma$ since it is technically a separate object. Note that $\gamma(0) = b$, because of the compatibility between $g$ and $f$. In this situation our lifting problem becomes that of finding a mapping $f_1 : P \times [0,1] \to E$ which is an extension of $f$ in the sense that $f_1(p, 0) = f(p)$ for all $p \in P$ and whose projection to the base is essentially $\gamma$ in the sense that $\pi(f_1(p, x)) = \gamma(x)$ for all $x \in [0,1]$ and all $p \in P$.

In other words, $f$ maps $P$ into a single fiber of $E$, and $f_1(\cdot, x)$ maps into the fiber in $E$ over $\gamma(x)$. As $x$ ranges through $[0,1]$ these fibers can move, but we can also return to the same fiber when $\gamma$ contains loops.

Let us explain how we can obtain how we can find such a lifting $f_1(p, x)$. If $E$ were just the product $S^1 \times T$ then we could pull $f$ along the parameter interval rigidly. Because of the twisting of our bundle we have to do something stronger, and it is convenient to go back to $\widetilde{E} = \mathbf{R} \times T$. Let $\widetilde{\gamma} : [0,1] \to \mathbf{R}$ be a continuous mapping which projects back to $\gamma$ under the canonical mapping from $\mathbf{R}$ to $\mathbf{R}/\mathbf{Z} = S^1$. This lifting $\widetilde{\gamma}$ of $\gamma$ is determined uniquely by its initial point $\widetilde{b} = \widetilde{\gamma}(0)$, which is a lifting of $b$ (but otherwise arbitrary). Our quotient mapping from $\widetilde{E}$ onto $E$ is a homeomorphism on each of the fibers, and so there is a mapping $\widetilde{f} : P \to \widetilde{E}$ which actually takes values in $\{\widetilde{b}\} \times T$ and which projects back down to $f$ when we project $\widetilde{E}$ onto $E$ using our quotient mapping.

In short, we can lift everything upstairs to $\widetilde{E}$. Since $\widetilde{E}$ is just a product we can define a mapping $\widetilde{f}_1 : P \times [0,1] \to \widetilde{E}$ in the obvious way, by setting

$$\widetilde{f}_1(p, x) = (\widetilde{\gamma}(x), \widetilde{\pi}_1(\widetilde{f}(p))),$$

where $\widetilde{\pi}_1 : \widetilde{E} \to T$ is the obvious projection. Thus up in $\widetilde{E}$ we are doing something quite trivial, we are simply sliding $\widetilde{f}$ along $\mathbf{R}$ rigidly.

Now we define $f_1 : P \times [0,1] \to E$ to simply be the composition of $\widetilde{f}_1 : P \times [0,1] \to \widetilde{E}$ with the quotient mapping from $\widetilde{E}$ onto $E$. It is easy to see that this choice of $f_1$ has the desired properties, namely that it is a continuous mapping which agrees with $f$ at the beginning and follows $\gamma$ in the base for the whole time interval $[0,1]$.

Now let us look at how $f_1(p, x)$ is distorted as $x$ runs from 0 to 1. Imagine that $\gamma$ moves at constant speed in $S^1$ but very fast, so that it wraps around $S^1$ many times. An integer number of times, $n$ times say, moving in the positive orientation and ending back at $b$. Then the lifting $\widetilde{\gamma}$ of $\gamma$ moves along at constant speed in $\mathbf{R}$, it starts at $\widetilde{b}$ and ends at $\widetilde{b} + n$.

So what does $f_1(\cdot, 1)$ look like when we go around the circle $n$ times? It looks like $f_1(\cdot, 0) = f$ acted on by $A^n$! That is, $f_1(\cdot, 1)$ and $f$ both map $P$ into the fiber $\pi^{-1}(b)$ in $E$, which is a copy of $T$. If we move the mappings back into $T$ so that we can look at them, then the transition from time 0 to time 1 is given by $A^n$. This is because each tour around $S^1$ corresponds to an application of $A$ on $T$. This follows from chasing definitions.

From our earlier analysis of $A$ (just before the second figure) we conclude that our mapping $f_1$ may undergo exponential stretching as we traverse the parameter interval $[0, 1]$. This is unavoidable and not simply an artifact of our construction. For instance, suppose that our initial mapping $f$ represents a loop in $T$ (i.e., $P$ is a circle). The ending mapping $f_1(\cdot, 1)$ represents another loop in $T$. No matter how we choose $f_1$ the homotopy class of our final loop in $T$ has to be the same as the one obtained from the construction above. Therefore the amount of winding that the final loop makes in terms of topology is simply determined by $A^n$, and can be exponentially large, as we have seen.

This finishes our concrete construction. Let us think about what it means. We are interested in feasibility for mappings into our space $E$. We might call such a mapping feasible if it is quite simple (e.g., if it does not wrap around too much), or if it can be obtained from a feasible mapping by a small perturbation. Thus feasibility is like being homotopic to something simple. One can then try to find complicated objects which are feasible with a short proof, and our construction of the lifting suggests an example of this.

Despite the difficulties of formalization this is certainly a situation worth considering. The problem of making liftings for fibrations is quite basic in topology. If one wants to try to come to terms with constructions in topology from the perspective of complexity of proofs it provides a good place to start.

# 9 Feasibility and continuous parameters

One might argue that the preceding example does not really mean anything, that the kinds of compression that we saw in earlier sections lose their sense in this continuous context. There is nothing wrong with that, it is basic issue to be understood. A point about our example is to see such issues in a context which is both fairly realistic in terms of ordinary mathematics and with enough structure for an idea like feasibility to be relevant.

A basic underlying problem, then, is this. Consider the idea of feasibility in an algebraic context, in terms of taking $B^n$ for large $n$ where $B$ is a matrix, say. In ordinary mathematics we are often allowed to pass to continuous constructions and take $B^t$ where $t$ is a real number. One can do this with matrices (positive definite matrices, for instance), and one is doing something very similar in solving an ordinary differential equation. Does anything remain in the idea of short proofs of feasibility in these situations?

One answer is simply "no". Consider the differential equation $y' = y$ which is solved by the exponential function. One can view this differential equation as a continuous version of recursion and the existence of its solution as a consequence of a "continuous" version of induction. (One can think of this as being analogous to proving by induction in arithmetic that exponential functions are defined everywhere.) For the notion of feasibility, it was important not to allow induction over $F$-formulas to avoid collapsing into triviality. A theory which provides the existence of exponential functions with continuous parameters, or solutions of differential equations, might arguably be too powerful to permit a meaningful notion of feasibility.

In ordinary mathematics one typically defines exponentials by summing infinite series, and one can find solutions of ordinary differential equations through approximation schemes which converge well. These constructions do not seem to be compatible with the short proofs of feasibility of numbers though. The short proofs try to operate from the "inside", speeding up the procedure by which we operate step by step, while the aforementioned con-structions of the exponential function and solutions of ordinary differential equations are more global in nature, dealing with a continuum of points at once.

One can try to make more "inner" definitions, such as

$$e^\alpha = \lim_{n \to \infty} \left(1 + \frac{\alpha}{n}\right)^n$$

for the exponential, and realizing solutions of ordinary differential equations as limits of discrete difference equations. For theoretical purposes these approaches are awkward, in large part because one must establish the existence of the limit. This is easier to do when one has the existence of the exponential or the solution to the differential equation established already by other means.

What about the topological example? In many situations one can treat it by solving an ordinary differential equation. On the other hand there is an extra discreteness to it. One does not need to know that the required family of mappings satisfies a differential equation but only that it exists and is continuous. Under reasonable conditions one does not even need a continuous family, one could replace the parameter interval $[0, 1]$ with a discrete set so long as the family that is constructed is "approximately" continuous, in the sense that each new step is a sufficiently small perturbation of the preceding one. One can then fill in the gaps automatically to get a continuous family parameterized by $[0, 1]$.

Roughly speaking one could say that for topology $(1 + \frac{\alpha}{n})^n$ is often practically as good as $e^\alpha$ when $n$ is large enough but still finite. This can allow the idea of feasibility to retain its relevance, perhaps in a more limited form. One has to be careful about what kind of operations are allowed.

Of course one is not prevented from making other notions of feasibility for approximating classical constructions in analysis (such as summing a series for the exponential function, or approximation schemes for finding solutions of ordinary differential equations). For this type of construction one typically deals with efficiency of approximation rather than exact values, as in the earlier and more algebraic examples. In topology a good approximation is often good enough already.

# 10 Remarks about compression and cuts

We have tried to explore some examples of *natural notions* of "feasibility", where there can be *short proofs*.

What does "short" really mean? In the examples there was some clear sense that the proofs were short compared to what one might expect, given the complexity of the particular object in question. But can we define "short" more abstractly, more invariantly, more objectively?

There is an obvious answer to this, in terms of the cut rule in sequent calculus. The reader who is unfamiliar with this may wish to consult [10] for an introduction to the cut rule and the combinatorics and complexity of cut elimination. The main point is that these short proofs all use cuts, and a natural way to measure the "shortness" is to ask how large a proof would have to be without the cut rule. We would then consider the *relative* sizes between proofs with and without cuts rather than the absolute sizes of proofs.

In many contexts in logic one can show that it is possible to eliminate cuts from a proof, but at great cost of expansion. See [18, 19, 21, 45]. There are examples known where the smallest proof without cuts is much larger than the smallest proof with cuts [46, 34, 35, 41, 42, 43, 27, 4]. The present discussion suggests that we view this phenomenon as a reflection of some kind of internal symmetry or structure. One can imagine results to the effect that there is much less to gain in the size of proofs using cuts in situations where the internal structure is not cooperative. In particular one might find that the gain is often not too dramatic. These issues are related to the "P = NP?" problem [14, 17, 28].

This may seem to be at odds with the normal experience of mathematicians, but then one has to face the issue that the situations that mathematicians can deal with are typically very special and enjoy a lot of internal structure.

Random objects cannot have short descriptions [30, 32, 11, 12]. This is well known and much-studied, although the basic point is clear by counting. The question then is whether certain types of descriptions are often much more efficient than others. In some contexts one knows that this is not the case, that the upper bounds on the best algorithmic description is about the same as the average. For formal proofs the general picture remains unclear.

*IHES*
*35 Route de Chartres*
*91440 Bures-sur-Yvette*
*France*

*A. Carbone*
*Institut für Algebra und Diskrete Mathematik*
*Technische Universität Wien*
*Wiedner Hauptstrasse 8-10/118*
*A-1040 Wien*
*Austria*

*S. Semmes*
*Department of Mathematics*
*Rice University*
*Houston Texas 77251*
*U.S.A.*